\documentclass{amsart}

\usepackage{xypic}
\input xy
\xyoption{all}
\usepackage{epsfig}
\usepackage{amsthm}
\usepackage{amssymb}
\usepackage{amsmath}
\usepackage{amscd}
\usepackage{color}
\usepackage[T1]{fontenc}

\newcommand{\bg}{\begin{equation}}
\newcommand{\ed}{\end{equation}}
\newcommand{\bga}{\begin{eqnarray}}
\newcommand{\eda}{\end{eqnarray}}
\newcommand{\pf}{\textbf{Proof:\ }}

\def\cbdu{\par{\raggedleft$\Box$\par}}

\newtheorem {Theorem}  {Theorem}

\numberwithin{Theorem}{section}

\newtheorem {Lemma}[Theorem]  {Lemma}

\theoremstyle{definition}
\newtheorem{Definition}[Theorem]{Definition}
\theoremstyle{remark}
\newtheorem{Remark}[Theorem]{\bf Remark}

\expandafter\chardef\csname pre amssym.def
at\endcsname=\the\catcode`\@ \catcode`\@=11
\def\undefine#1{\let#1\undefined}
\def\newsymbol#1#2#3#4#5{\let\next@\relax
 \ifnum#2=\@ne\let\next@\msafam@\else
 \ifnum#2=\tw@\let\next@\msbfam@\fi\fi
 \mathchardef#1="#3\next@#4#5}
\def\mathhexbox@#1#2#3{\relax
 \ifmmode\mathpalette{}{\m@th\mathchar"#1#2#3}%
 \else\leavevmode\hbox{$\m@th\mathchar"#1#2#3$}\fi}
\def\hexnumber@#1{\ifcase#1 0\or 1\or 2\or 3\or 4\or 5\or 6\or 7\or 8\or
 9\or A\or B\or C\or D\or E\or F\fi}

\font\teneufm=eufm10 \font\seveneufm=eufm7 \font\fiveeufm=eufm5
\newfam\eufmfam
\textfont\eufmfam=\teneufm \scriptfont\eufmfam=\seveneufm
\scriptscriptfont\eufmfam=\fiveeufm

\catcode`\@=\csname pre amssym.def at\endcsname

\newcounter{remark}
\newcommand{\supp}{{\mathit supp}\,}

\newcommand{\e}{\epsilon}

\renewcommand{\th}{\theta}

\newcommand{\R}{\mathbf{R}}

\renewcommand{\div}{\mbox{div}}

\def  \R   {{\mathbb R}}

\def  \12  {{\frac{1}{2}}}

\def\build#1_#2^#3{\mathrel{\mathop{\kern 0pt#1}\limits_{#2}^{#3}}}

\begin{document}
%\currannalsline{0}{2006}

\title[Regularity criterion and energy conservation]{Regularity criterion and energy conservation for the supercritical Quasi-Geostrophic equation}

%\author{hello}

%\author [Alexey Cheskidov]{Alexey Cheskidov}
%\address{Department of Mathematics, Stat. and Comp.Sci.,  University of Illinois Chicago, Chicago, IL 60607,USA}
%\email{acheskid@math.uic.edu} 
\author [Mimi Dai]{Mimi Dai}
\address{Department of Mathematics, University of Illinois at Chicago, IL 60607,USA}
\email{mdai@uic.edu} 

%\thanks{The work of Alexey Cheskidov was partially supported by NSF Grant DMS--1108864.}

%%%use \Proof instead of \begin{proof}

%%%% use \Endproof instead of \end{proof}

%%%% use \references {999} instead of \begin{thebibliography}{99}

%%%%used \Endrefs instead of \end{thebibliography}

\begin{abstract}
This paper studies the regularity and energy conservation problems for the 2D supercritical quasi-geostrophic (SQG) equation. We apply an approach of splitting the dissipation wavenumber %combined with an estimate of the energy flux 
to obtain a new regularity condition which is weaker than all the Prodi-Serrin type regularity conditions. Moreover, we prove that any viscosity solution of the supercritical SQG in $L^2(0,T; B^{1/2}_{2,c(\mathbb N)})$ satisfies energy equality.
%This splitting approach was originally introduced by Cheskidov and Shvydkoy \cite{CSr} to establish optimal regularity criteria for both the 3D Navier-Stokes and Euler equations. 
%Onsagar's conjecture for supercritical SQG;
%Small data result for supercritical SQG in the largest critical space $\dot B^{1-\alpha}_{\infty,\infty}$.

\bigskip


\hspace{0.02cm}CLASSIFICATION CODE: 76D03, 35Q35.
\end{abstract}

\maketitle

\section{Introduction}

In this paper we consider the two dimensional supercritical quasi-geostrophic (SQG) equation

\begin{equation}\begin{split}\label{QG}
\theta_t+u\cdot\nabla \theta+\kappa\Lambda^\alpha\theta =0,\\
u=R^\perp\theta,
\end{split}
\end{equation}
where $0< \alpha\leq 2$, $\kappa>0$, $\Lambda=\sqrt{-\Delta}$ is the Zygmund operator, and
\bg\notag
R^\perp\theta=\Lambda^{-1}(-\partial_2\theta,\partial_1\theta).
\ed
The equations are considered in $\R^2\times(0, T)$ with $T>\infty$. The scalar function $\theta$ represents the potential temperature and the vector function $u$ represents the fluid velocity. Notice that if $\theta(x,t)$ solves (\ref{QG}) with initial data $\theta_0(x)$, the function $\theta_\lambda=\lambda^{\alpha-1}\theta(\lambda x, \lambda^\alpha t)$ also solves (\ref{QG}) with initial data $\theta_0(\lambda x)$. The space which is invariant under such scaling is called a critical space. For example, $\dot H^{2-\alpha}$ and $\dot B_{\infty,\infty}^{1-\alpha}$ are critical, and the latter one is the largest critical space of (\ref{QG}).

Equation (\ref{QG}) with $\alpha=1$ %corresponds to a fundamental model of the quasi-geostrophic equation. It 
describes the evolution of the surface temperature field in a rapidly rotating and stably stratified fluid with potential velocity. As pointed out in \cite{CMT}, this equation attracts interest of  scientists and mathematicians due to two major reasons: it is a fundamental model for the actual geophysical flows with applications in atmosphere and oceanography study; from the mathematical point of view, the behavior of strongly nonlinear solutions to (\ref{QG}) with $\kappa=0$ in 2D and the behavior of potentially singular solutions to the Euler's equation in 3D are strikingly analogous which has been justified both analytically and numerically. For literature the readers are refereed to \cite{CMT, CW, Pe} and the references therein. 

Equation (\ref{QG}) is usually referred as supercritical, critical and subcritical SQG for $0<\alpha<1$, $\alpha=1$ and $1<\alpha\leq2$ respectively, although it is an open problem whether a dramatic change in the behavior of solutions occurs for the case of dissipation power less than $1$. The global regularity problem of the critical SQG equation has been very challenging due to the balance of the nonlinear term and the dissipative term
in (\ref{QG}). This problem is resolved now by Kieslev, Nazarov and Volberg \cite{KNV}, Caffarelli and Vasseur \cite{CV}, Kieslev and Nazarov \cite{KN09} and Constantin and Vicol \cite{CVicol} independently, using different sophisticated methods. 
%The main ingredient of the proof of global regularity in \cite{KN} relies on constructing a special family of Lipschitz moduli of continuity that are preserved by the dissipative evolution. The proof of \cite{CaV} applies the ideas of De Giorgi iteration method to the nonlocal parabolic equation. As a first step, the authors make use of the interplay between $|\Lambda ^{\frac 12}\theta|$ and $|\theta|$ to prove that a weak solution in $L^2$ is bounded in $L^\infty$. The second step, with a more delicate analysis, such a weak solution with finite $L^\infty$ norm is proved to be H\"older continuous. Later, to find a bridge between the proofs of \cite{KN} and \cite{CaV}, Kieslev and Nazarov \cite{KN09} reproved the global regularity using a completely different method which applied elementary tools to control the H\"older norms by choosing a suitable family of test functions (Hardy molecules). Very recently, Constantin and Vicol \cite{CV} proposed another new proof of the global regularity which provides a more transparent way to see that the dissipation is dominating the nonlinear term in the critical SQG equation. The main tool is a nonlinear maximum principle which introduces nonlinear lower bounds for the linear nonlocal operator $\Lambda$.

The global regularity problem of the supercritical SQG equation (\ref{QG}) ($0<\alpha<1$) remains open. In \cite{CMT}, Constantin, Majda and Tabak established the first regularity criterion when $\kappa=0$:
\[\lim\sup_{t\to T}\|\theta(t)\|_{H^m}<\infty \quad \mbox { if and only if } \quad \int_0^T\|\nabla^\perp\theta(t)\|_{L^\infty}\, dt<\infty\]
for $m>2$. A Prodi-Serrin type regularity criterion was obtained by Chae \cite{Ch06}: if the solution $\theta(x,t)$ satisfies 
\[\theta\in L^r(0,T; L^p(\mathbb R^2)), \qquad \mbox { for } \quad \frac2p+\frac\alpha r\leq\alpha, \quad\frac2\alpha<p<\infty\]
then there is no singularity up to time $T$.
Many works have been devoted to extend and improve the above regularity criteria, for instance, see \cite{Ch03, CW09, CW08, DP, FGN, Xi, Ya, Yu, Yuan, Yuan8, ZV}. In particular, in \cite{CW08} the authors proved that a weak solution in the H\"older space $C^\delta$ with $\delta>1-\alpha$ is actually classical, see Theorem \ref{CW}. The result of \cite{Ch03} improved the Constantin-Majda-Tabak criterion by using Triebel-Lizorkin spaces. In \cite{DP}, the authors obtained the regularity for a weak solution $\theta$ which belongs to $\theta\in L^s((0,T); B^\gamma_{p,\infty})$ with $\gamma=\frac2p+1-\alpha+\frac\alpha s$, for $1\leq s<\infty$ and $2< p<\infty$, which is a Prodi-Serrin type regularity criterion. Besides, the results of \cite{FGN, Xi, Yuan, Yuan8} can be also considered as extensions of the Prodi-Serrin type regularity criteria. The work of \cite{Ya} impose regularity condition only on partial derivative of the solutions. While the authors of \cite{Yu, ZV} obtained global regularity by imposing conditions on the initial data. Moreover, in this area,
eventual regularity results for a slightly supercritical SQG and supercritical SQG are obtained in \cite{Sil, Dab}, respectively. And later the global regularity for a slightly supercritical active scalar equation is obtained in \cite{DKLV}.

Among the conditional regularity results mentioned above, the work of \cite{CW08, DP} are known to be sharp in the case of linear drift-diffusion equations (see \cite{SVZ}). In the first part of the paper, we will prove a new regularity criterion for the supercritical SQG by showing that under such new condition a viscosity solution is in the H\"older space $C^\delta$ with $\delta>1-\alpha$, and hence classical (see \cite{CW08}). We will also show that our regularity condition is weaker than all the Prodi-Serrin type regularity conditions, particularly including the one of \cite{DP}. 
Our main result states as follows.

\begin{Theorem}\label{thm}
Let $0<\alpha<1$. Let $\theta$ be a viscosity solution to (\ref{QG}) on $[0,T]$. Assume that $\theta(t)$ is regular on $[0,T)$, and 
\bg\notag%\label{criteria}
\int_0^T\|\nabla\theta_{\leq Q(t)}(t)\|_{B^0_{\infty,\infty}}dt<\infty, \qquad \mbox { for a certain
 number} \qquad Q(t).
\ed
Then $\theta(t)$ is regular on $[0,T]$.
\end{Theorem}

\begin{Remark} In fact the statement holds true for any weak solution satisfying the truncated energy estimates as in \cite{CV}.
\end{Remark}

The number $Q(t)$ will be defined in Section \ref{sec:reg}. Roughly speaking, this criterion says, if the solution at low modes is bounded in $L^1(0,T; B^1_{\infty,\infty})$, then it is regular up to time $T$. The intuition is that the linear dissipation term $\Lambda^\alpha\theta$ dominates at high modes. Following the idea of the work \cite{CSr} for the NSE and Euler equation, we split the dissipation by an appropriate wavenumber $\lambda_{Q(t)}$ (see notation in Section \ref{sec:pre}). We show that above $\lambda_{Q(t)}$, the nonlinear interaction is dominated by the linear term $\Lambda^\alpha\theta$; while below $\lambda_{Q(t)}$, the nonlinear interaction is controlled due to the assumption on the low modes as in Theorem \ref{thm}. 

The second part of the paper concerns the energy conservation problem for the supercritical SQG (\ref{QG}). The well known Onsager's conjecture (see \cite{On}) addresses the energy conservation for the Euler equation which states: any weak solution of the Euler equation with H\"older continuity $s>1/3$ conserves energy and dissipates energy otherwise. In \cite{CCFS}, the authors proved that the energy is conserved for any weak solution in the Besov space $B^{1/3}_{3, c(\mathbb N)}$ (slightly smaller than $B^{1/3}_{3, \infty}$) in which the ``H\"older exponent" is exactly $1/3$. It is also shown that the space $B^{1/3}_{3, c(\mathbb N)}$ is sharp in the context of no anomalous dissipation. Compared to the Euler equation, the SQG has a nonlinear term $R^\perp\theta\cdot\nabla\theta$ with the same degree of derivative. One may expect that an analogous result holds for the SQG: any visicosity solution of the supercritical SQG in $B^{1/3}_{3, c(\mathbb N)}$ satisfies energy equality. However, it is shown in \cite{CV} that any viscosity solution of the critical SQG is in $L^\infty([t_0,\infty);L^\infty(\mathbb R^2))$ for every $t_0>0$. This result holds for the supercritical SQG (even with an external force) as well (see proofs in \cite{CD, CDatt}). Thanks to this fact, we are able to prove that any viscosity solution $\theta$ of the supercritical SQG satisfies energy equality provided $\theta\in L^2(0,T; B^{1/2}_{2, c(\mathbb N)})$. Notice that $L^2(0,T; B^{1/3}_{3, c(\mathbb N)})\subset  L^2(0,T; B^{1/2}_{2, c(\mathbb N)})$. The main result states as follows.

\begin{Theorem}\label{thm-energy}
Let $\theta\in C_w([0,T]; L^2)$ be a viscosity solution to the supercritical SQG (\ref{QG}) with $0<\alpha<1$. If additionally, we assume $\theta\in L^2(0,T; B^{1/2}_{2, c(\mathbb N)})$, then the solution $\theta$ satisfies energy equality.
\end{Theorem}

The rest of the paper is organized as follows: in Section \ref{sec:pre} we introduce some notations, recall the Littlewood-Paley decomposition theory briefly, and recall a regularity result for the supercritical SQG; Section \ref{sec:reg} and \ref{sec-on} are devoted to proving Theorem \ref{thm} and \ref{thm-energy}, respectively.

\bigskip

\section{Preliminaries}
\label{sec:pre}

\subsection{Notation}
\label{sec:notation}
We denote by $A\lesssim B$ an estimate of the form $A\leq C B$ with
some absolute constant $C$, and by $A\sim B$ an estimate of the form $C_1
B\leq A\leq C_2 B$ with some absolute constants $C_1$, $C_2$. We agree that
$\|\cdot\|_p=\|\cdot\|_{L^p}$. While $(\cdot, \cdot)$ stands for the $L^2$-inner product.

\subsection{Littlewood-Paley decomposition}
\label{sec:LPD}
The techniques presented in this paper rely strongly on the Littlewood-Paley decomposition. Thus we recall the Littlewood-Paley decomposition theory here briefly. For a more detailed description on this theory we refer the readers to the books by Bahouri, Chemin and Danchin \cite{BCD} and Grafakos \cite{Gr}. 

In this subsection, $u$ denotes a general function, which should not be considered as the $u$ in (\ref{QG}).

Denote $\lambda_q=2^q$ for integers $q$. We choose a radial function $\chi\in C_0^\infty(\R^n)$ as
\begin{equation}\notag
\chi(\xi)=
\begin{cases}
1, \ \ \mbox { for } |\xi|\leq\frac{3}{4}\\
0, \ \ \mbox { for } |\xi|\geq 1.
\end{cases}
\end{equation}
Let 
\begin{equation}\notag
\varphi(\xi)=\chi(\frac{\xi}{2})-\chi(\xi),  \qquad \mbox { and } \qquad
\varphi_q(\xi)=
\begin{cases}
\varphi(\lambda_q^{-1}\xi)  \ \ \ \mbox { for } q\geq 0,\\
\chi(\xi) \ \ \ \mbox { for } q=-1.
\end{cases}
\end{equation}
For a tempered distribution vector field $u$ we define a Littlewood-Paley decomposition
\begin{equation}\notag%\label{eq:LPP}
\begin{cases}
&h=\mathcal F^{-1}\varphi, \qquad \tilde h=\mathcal F^{-1}\chi,\\
&u_q:=\Delta_qu=\mathcal F^{-1}(\varphi(\lambda_q^{-1}\xi)\mathcal Fu)=\lambda_q^n\int h(\lambda_qy)u(x-y)dy,  \qquad \mbox { for }  q\geq 0,\\
& u_{-1}=\mathcal F^{-1}(\chi(\xi)\mathcal Fu)=\int \tilde h(y)u(x-y)dy.
\end{cases}
\end{equation}
Recall that the following identity
\bg\notag
u=\sum_{q=-1}^\infty u_q
\ed
holds in the distribution sense. Essentially the sequence of the smooth functions $\varphi_q$ forms a dyadic partition of the unit. 

To simplify the notation, we denote
\bg\notag
u_{\leq Q}=\sum_{q=-1}^Qu_q, \qquad \tilde u_q=u_{q-1}+u_q+u_{q+1}.
\ed
By the definition of $\varphi_q$, we have
$$\supp (\hat u_q)\cap\supp (\hat u_p)=\emptyset  \qquad \mbox { if }  |p-q|\geq 2.$$

By the Littlewood-Paley decomposition we define the inhomogeneous Besov spaces $B_{p,r}^{s}$.
\begin{Definition}
Let $s\in \mathbb R$, and $1\leq p, r\leq \infty$. Then
$$
\|u\|_{B_{p,r}^{s}}=\|u_q\|_p+\left\|\left(\lambda_q^s\|u_q\|_p\right)_{q\in\mathbb N}\right\|_{l^r(\mathbb N)}
$$
is the inhomogeneous Besov norm. The inhomogeneous Besov space $B_{p,r}^{s}$ is the space of tempered distributions $u$ such that the norm $\|u\|_{B_{p,r}^{s}}$ is finite.
\end{Definition}
Specially, a tempered distribution $u$ belongs to $B_{p, \infty}^{s}$ if and only if
$$
\|u\|_{B_{p, \infty}^{s}}=\sup_{q}\lambda_q^s\|u_q\|_p<\infty.
$$

\begin{Definition}
We define $B^{1/2}_{2,c(\mathbb N)}$ as the class of all tempered distributions $u$ for which
\[\lim_{q\to\infty}\lambda_q^{1/2}\|u_q\|_2=0.\]
The space $B^{1/2}_{2,c(\mathbb N)}$ is endowed with the norm inherited from $B^{1/2}_{2,\infty}$.
\end{Definition}

We recall two inequalities for the dyadic blocks of the Littlewood-Paley decomposition in the following.
\begin{Lemma}\label{le:bern}(Bernstein's inequality) \cite{L}
Let $n$ be the space dimension and $r\geq s\geq 1$. Then for all tempered distributions $u$, 
\bg\notag%\label{Bern}
\|u_q\|_{r}\lesssim \lambda_q^{n(\frac{1}{s}-\frac{1}{r})}\|u_q\|_{s}.
\ed
\end{Lemma}

\begin{Lemma}\label{lp}
Assume $2<l<\infty$ and $0\leq \alpha\leq 2$. Then
\begin{equation}\notag
l\int_{\mathbb R^n}u_q\Lambda^\alpha u_q|u_q|^{l-2} \, dx\gtrsim  \lambda_q^\alpha\|u_q\|_l^l.
\end{equation}
%for some constant $C$ depending on $n$ and $\alpha$.
\end{Lemma}
For a proof of Lemma \ref{lp}, see \cite{CC, CMZ, Wu}.

\bigskip

\subsection{Weak solution, viscosity solution and H\"older regularity}
\label{sec:sol}
We recall the standard definition of weak solutions, viscosity solutions, and a regularity result for the supercritical SQG.

\begin{Definition}\label{def:weak}
A Leray-Hopf weak solution of (\ref{QG}) on $[0,T]$ (or $[0, \infty)$ if $T=\infty$) is a function $\theta\in C_w([0,T]; L^2(\mathbb R^2))$
with $\theta(x,0)=\theta_0$, satisfying, for all test functions $\phi\in C_0^\infty([0,T]\times\mathbb R^2)$ with $\nabla_x\cdot \phi=0$
\begin{equation}\notag%\label{eq:uweak}
\begin{split}
&\int_0^t(\theta(s), \partial_s\phi(s))+\kappa(\Lambda^{\alpha/2}\theta(s), \Lambda^{\alpha/2}\phi(s))+ (\theta(s)\cdot\nabla\phi(s), \theta(s))ds\\
=&(\theta(t), \phi(t))-(\theta_0, \phi(0))
\end{split}
\end{equation}
and moreover, the following energy inequality 
\begin{equation}\notag%\label{energy}
\|\theta(t)\|_2^2+2\kappa\int_{t_0}^t\|\nabla \theta(s)\|_2^2ds
\leq \|u(t_0)\|_2^2
\end{equation}
is satisfied for almost all $t_0\in(0, T)$ and all $t\in(t_0, T]$.
\end{Definition}

\begin{Definition}
A weak solution $\th(t)$ on $[0,T]$ is called a viscosity solution if there exist sequences $\e_n \to 0$ and $\th_n(t)$ satisfying
\begin{equation}\begin{split}\notag%\label{VQG}
\frac{\partial\theta_n}{\partial t}+u_n\cdot\nabla \theta_n+\kappa\Lambda\theta_n + \e_n \Delta \th_n=0,\\
u_n=R^\perp\theta_n,
%\theta(x,0)=\theta_0,
\end{split}
\end{equation}
such that $\th_n \to \theta$ in $C_\mathrm{w}([0,T];L^2)$. 
\end{Definition}

Standard arguments imply that for any initial data $\th_0 \in L^2$ there exists a viscosity solution $\th(t)$ of \eqref{QG} on $[0,\infty)$ with $\th(0)=\th_0$ (see \cite{CC}, for example).

Constantin and Wu \cite{CW08} proved the following regularity result for the supercritical SQG (\ref{QG}).

\begin{Theorem}\label{CW}
Let $\theta$ be a Leray-Hopf weak solution of (\ref{QG}) with $0<\alpha<1$. Let $\delta>1-\alpha$ and  $0<t_0<t<\infty$. If 
\[\theta\in L^\infty([t_0,t]; C^\delta(\mathbb R^2)),\]
then 
\[\theta\in C^\infty((t_0,t]\times\mathbb R^2).\]
\end{Theorem}

\bigskip

\section{Regularity criterion}
\label{sec:reg}

In this section we prove the regularity criterion stated in Theorem \ref{thm}, and show that this regularity condition is weaker than all the Prodi-Serrin type regularity conditions.

{\bf{Proof of Theorem \ref{thm}:}}
 Let $\theta(t)$ be a weak solution of (\ref{QG}) on $[0,T]$. We adopt the notations and idea from \cite{CSr} and define the dissipation wavenumber as
\begin{equation}\label{wave}
\Lambda(t)=\min \left\{\lambda_q:\lambda_p^{1-\alpha}\|\theta_p(t)\|_\infty<c_0\kappa, \quad \forall p>q\geq 1 \right\},  
\end{equation}
where $c_0$ is an absolute constant which will be determined later. Let $Q(t)\in\mathbb N$ be such that $\lambda_{Q(t)}=\Lambda(t)$.  It follows immediately that
\bg\notag%\label{ineq:wave}
\|\theta_{Q(t)}(t)\|_\infty\geq c_0\kappa\Lambda(t)^{-1+\alpha},
\ed
provided $1<\Lambda(t)<\infty$. We consider the function 
\[
f(t)=\|\theta_{\leq Q(t)}(t)\|_{B^1_{\infty,\infty}}=\sup_{q\leq Q(t)}\lambda_q\|\theta_q(t)\|_\infty.
\]
%which is equivalent to the form in term of vorticity $\omega=\nabla\times u$,
%$$
%f(t)\sim \|\omega_{\leq Q(t)}(t)\|_{B^0_{\infty,\infty}}.
%$$
The idea is to prove that $\|\theta(t)\|_{B^s_{l,l}}$ is uniformly bounded on $[0,T)$ for some large integer $l$ and $s\in(0,1)$ provided $f\in L^1(0,T)$. Notice that for $0<s<1$ and $sl>2$,  we have the embedding $B^s_{l,l}\subset C^{0, s-\frac2l}$, see \cite{DPV}. Choose  large enough $l$ such that 
\begin{equation}\label{para1}
s-\frac2l>1-\alpha,
\end{equation}  
it then follows from Theorem \ref{CW} that, $\theta$ is regular on $[0,T]$.  In addition, we assume
\begin{equation}\label{para2}
1-s-\frac\alpha l>0
\end{equation}
to carry through the estimate for $\|\theta(t)\|_{B^s_{l,l}}$. Notice that the two conditions (\ref{para1}) and (\ref{para2}) are compatible. Indeed, for any $\alpha\in(0,1)$, one can choose large enough $l$ such that $\frac\alpha l<1-s<\alpha-\frac2l$.% which implies (\ref{para1}) and (\ref{para2}) are compatible.
%On the other hand, it is clear that $W^{s,l}\subset \dot H^{2-\alpha}$ which is a critical space for (\ref{QG}). 

Now we prove $\|\theta(t)\|_{B^s_{l,l}}$ is uniformly bounded on $[0,T)$. Since $\theta(t)$ is regular on $(0,T)$, projecting (\ref{QG}) onto the $q-$th shell, testing it with $l\lambda_q^{sl}\theta_q|\theta_q|^{l-2}$, summing over $q\geq -1$, and applying Lemma \ref{lp} yields 
\begin{equation}\label{ineq-q}
\frac{d}{dt}\sum_{q\geq -1}\lambda_q^{sl}\|\theta_q\|_l^l\leq -C\kappa \sum_{q\geq-1}\lambda_q^{sl+\alpha}\|\theta_q\|_l^l+l\sum_{q\geq -1}\lambda_q^{sl}\int_{\R^3}\Delta_q(u\cdot\nabla \theta)\theta_q|\theta_q|^{l-2}\, dx.
\end{equation}
Using Bony's notation ofoduct, we write
\begin{equation}\notag
\begin{split}
\Delta_q(u\cdot\nabla \theta)=&\sum_{|q-p|\leq 2}\Delta_q(u_{\leq{p-2}}\cdot\nabla \theta_p)+
\sum_{|q-p|\leq 2}\Delta_q(u_{p}\cdot\nabla \theta_{\leq{p-2}})\\
&+\sum_{p\geq q-2}\Delta_q(u_p\cdot\nabla\tilde \theta_p).
\end{split}
\end{equation}
Recall the commutator notation
\begin{equation}\notag
[\Delta_q, u_{\leq{p-2}}\cdot\nabla]\theta_p=\Delta_q(u_{\leq{p-2}}\cdot\nabla \theta_p)-u_{\leq{p-2}}\cdot\nabla \Delta_q\theta_p.
\end{equation}
Thus, we decompose the integral
\begin{equation}\label{eq-i}
\begin{split}
&l\sum_{q\geq -1}\lambda_q^{sl}\int_{\R^3}\Delta_q(u\cdot\nabla \theta)\theta_q|\theta_q|^{l-2}\, dx\\
=&l\sum_{q\geq -1}\sum_{|q-p|\leq 2}\lambda_q^{sl}\int_{\R^3}\Delta_q(u_{\leq{p-2}}\cdot\nabla\theta_p) \theta_q|\theta_q|^{l-2}\, dx\\
&+l\sum_{q\geq -1}\sum_{|q-p|\leq 2}\lambda_q^{sl}\int_{\R^3}\Delta_q(u_{p}\cdot\nabla \theta_{\leq{p-2}})\theta_q|\theta_q|^{l-2}\, dx\\
&+l\sum_{q\geq -1}\sum_{p\geq q-2}\lambda_q^{sl}\int_{\R^3}\Delta_q(u_p\cdot\nabla \tilde\theta_p)\theta_q|\theta_q|^{l-2}\, dx\\
=&I_{1}+I_{2}+I_{3}.
\end{split}
\end{equation}
Using the commutator notation, $I_1$ can be further decomposed as
\begin{equation}\notag
\begin{split}
I_1=&l\sum_{q\geq -1}\sum_{|q-p|\leq 2}\lambda_q^{sl}\int_{\R^3}[\Delta_q,u_{\leq{p-2}}\cdot\nabla] \theta_p\theta_q|\theta_q|^{l-2}\, dx\\
&+l\sum_{q\geq -1}\lambda_q^{sl}\int_{\R^3}u_{\leq q-2}\cdot\nabla \theta_q \theta_q|\theta_q|^{l-2}\, dx\\
&+l\sum_{q\geq -1}\sum_{|q-p|\leq 2}\lambda_q^{sl}\int_{\R^3}(u_{\leq{p-2}}-u_{\leq q-2})\cdot\nabla\Delta_q\theta_p \theta_q|\theta_q|^{l-2}\, dx\\
=&I_{11}+I_{12}+I_{13},
\end{split}
\end{equation}
where we used $\sum_{|q-p|\leq 2}\Delta_q \theta_p=\theta_q$.
One can see that $I_{12}=0$, using integration by parts and the fact $\div\, u_{\leq q-2}=0$. By the definition of $\Delta_q$,
\begin{equation}\notag
\begin{split}
[\Delta_q,u_{\leq{p-2}}\cdot\nabla] \theta_q=&\lambda_q^3\int_{\R^3}h(\lambda_q(x-y))\left(u_{\leq p-2}(y)-u_{\leq p-2}(x)\right)\nabla \theta_q(y)\,dy\\
=&-\lambda_q^3\int_{\R^3}\nabla h(\lambda_q(x-y))\left(u_{\leq p-2}(y)-u_{\leq p-2}(x)\right) \theta_q(y)\,dy.
\end{split}
\end{equation}
By Young's inequality, for $1\leq r\leq \infty$,
\begin{equation}\notag
\begin{split}
&\|[\Delta_q,u_{\leq{p-2}}\cdot\nabla] \theta_q\|_{r}\\
\lesssim &\|\nabla u_{\leq p-2}\|_{r}\|\theta_q\|_{\infty}\left|\lambda_q^3\int_{\R^3}|x-y|\nabla h(\lambda_q(x-y))\, dy\right|\\
\lesssim &\|\nabla u_{\leq p-2}\|_{r}\|\theta_q\|_{\infty}.
\end{split}
\end{equation}
Thus, by taking $p=q$, $I_{11}$ can be estimated as, 
\begin{equation}\notag
\begin{split}
|I_{11}|\lesssim &l\sum_{q\geq-1}\lambda_q^{sl}\|[\Delta_q,u_{\leq{q-2}}\cdot\nabla] \theta_q\|_l\|\theta_q\|_l^{l-1}\\
\lesssim&l\sum_{q> Q}\lambda_q^{sl}\|\theta_q\|_\infty\|\theta_q\|_l^{l-1}\sum_{p\leq q-2}\lambda_p\|u_p\|_l\\
&+l\sum_{-1\leq q\leq Q}\lambda_q^{sl}\|\theta_q\|_\infty\|\theta_q\|_l^{l-1}\sum_{p\leq q-2}\lambda_p\|u_p\|_l\\
\equiv& A+B,
\end{split}
\end{equation}
while using the fact $\|u\|_p\lesssim \|\theta\|_p$, it follows
\begin{equation}\notag
\begin{split}
A\lesssim&c_0\kappa l\sum_{q> Q}\lambda_q^{sl-1+\alpha}\|\theta_q\|_l^{l-1}\sum_{p\leq q-2}\lambda_p\|u_p\|_l\\
\lesssim&c_0\kappa l\sum_{q> Q}\lambda_q^{sl+\alpha-\frac {sl+\alpha}{l}}\|\theta_q\|_l^{l-1}\sum_{p\leq q-2}\lambda_p^{\frac{sl+\alpha}l}\|u_p\|_l\lambda_{q-p}^{-1+\frac{sl+\alpha}{l}}\\
\lesssim&c_0\kappa l\sum_{q> Q}\lambda_q^{sl+\alpha}\|\theta_q\|_l^{l}+c_0\kappa l\sum_{q> Q}\left(\sum_{p\leq q-2}\lambda_p^{\frac{sl+\alpha}l}\|u_p\|_l\lambda_{q-p}^{-1+\frac{sl+\alpha}{l}}\right)^l\\
\lesssim&c_0\kappa l\sum_{q\geq -1}\lambda_q^{sl+\alpha}\|\theta_q\|_l^{l}
\end{split}
\end{equation}
due to the fact $1-s-\frac{\alpha}{l}>0$ in assumption (\ref{para2}); and 
\begin{equation}\notag
\begin{split}
B\lesssim &lf(t)\sum_{-1\leq q\leq Q}\lambda_q^{sl-1}\|\theta_q\|_l^{l-1}\sum_{p\leq q-2}\lambda_p\|u_p\|_l\\
\lesssim &lf(t)\sum_{-1\leq q\leq Q}\lambda_q^{sl}\|\theta_q\|_l^{l-1}\sum_{p\leq q-2}\lambda_{p-q}\|u_p\|_l\\
\lesssim &lf(t)\sum_{-1\leq q\leq Q}\lambda_q^{sl}\|\theta_q\|_l^{l}.
\end{split}
\end{equation}
The term $I_{13}$ is estimated as, 
\begin{equation}\notag
\begin{split}
|I_{13}|\lesssim &l\sum_{q\geq-1}\lambda_q^{sl}\int_{\mathbb R^3}|u_q| |\nabla\theta_q| |\theta_q|^{l-1}\, dx\\
\lesssim&l\sum_{q> Q}\lambda_q^{sl+1}\|u_q\|_l\|\theta_q\|_\infty\|\theta_q\|_l^{l-1}
+l\sum_{-1\leq q\leq Q}\lambda_q^{sl+1}\|u_q\|_l\|\theta_q\|_\infty\|\theta_q\|_l^{l-1}\\
\lesssim &c_0\kappa l\sum_{q> Q}\lambda_q^{sl+1-1+\alpha}\|u_q\|_l\|\theta_q\|_l^{l-1}
+lf(t)\sum_{-1\leq q\leq Q}\lambda_q^{sl}\|u_q\|_l\|\theta_q\|_l^{l-1}\\
\lesssim &c_0\kappa l\sum_{q> Q}\lambda_q^{sl+\alpha}\|\theta_q\|_l^l+lf(t)\sum_{-1\leq q\leq Q}\lambda_q^{sl}\|\theta_q\|_l^l.
\end{split}
\end{equation}
Therefore, we have
\begin{equation}\label{est-i1}
|I_1|\lesssim c_0\kappa l\sum_{q\geq -1}\lambda_q^{sl+\alpha}\|\theta_q\|_l^{l}+lf(t)\sum_{-1\leq q\leq Q}\lambda_q^{sl}\|\theta_q\|_l^{l}.
\end{equation}
While for $I_2$, 
\begin{equation}\notag
\begin{split}
|I_2|\lesssim& l\sum_{q> Q}\lambda_q^{sl}\|u_q\|_\infty\|\nabla \theta_{\leq q-1}\|_l\|\theta_q\|_l^{l-1}
+l\sum_{-1\leq q\leq Q}\lambda_q^{sl}\|u_q\|_\infty\|\nabla \theta_{\leq q-1}\|_l\|\theta_q\|_l^{l-1}.
\end{split}
\end{equation}
Notice that $\|u_q\|_\infty\lesssim \|\theta_q\|_\infty$ for any $q$, we conclude that $I_2$ has the same estimate as $I_{11}$. Thus, %for $1+\frac2r-s-\frac{\alpha}{l}>0$,
\begin{equation}\label{est-i2}
|I_2|\lesssim c_0\kappa l\sum_{q\geq -1}\lambda_q^{sl+\alpha}\|\theta_q\|_l^{l}+lf(t)\sum_{-1\leq q\leq Q}\lambda_q^{sl}\|\theta_q\|_l^{l}.
\end{equation}
While $I_3$ is estimated as
\begin{equation}\notag
\begin{split}
|I_3|\lesssim&\left|l^2\sum_{q\geq -1}\lambda_q^{sl}\sum_{p\geq q}\int_{\mathbb R^3}\Delta_q(u_p \theta_p)\nabla \theta_q \theta_q^{l-2}\, dx\right|\\
\lesssim &l^2\sum_{p\geq-1}\sum_{-1\leq q\leq p}\lambda_q^{sl}\int_{\mathbb R^3}|\Delta_q(u_p\theta_p)||\nabla\theta_q||\theta_q|^{l-2}\, dx\\
\lesssim & l^2\sum_{p> Q}\|u_p\|_l\|\theta_p\|_\infty\sum_{-1\leq q\leq p}\lambda_q^{sl+1}\|\theta_q\|_l^{l-1}\\
&+l^2\sum_{-1\leq p\leq Q}\|u_p\|_l\|\theta_p\|_\infty\sum_{-1\leq q\leq p}\lambda_q^{sl+1}\|\theta_q\|_l^{l-1}\\
\equiv &C+D.
\end{split}
\end{equation}
Similar analysis gives that
\begin{equation}\notag
\begin{split}
C\lesssim &c_0\kappa l^2\sum_{p> Q}\lambda_p^{-1+\alpha}\|u_p\|_l
\sum_{-1\leq q\leq p}\lambda_q^{sl+1}\|\theta_q\|_l^{l-1}\\
\lesssim &c_0\kappa l^2\sum_{p> Q}\lambda_p^{\frac{sl+\alpha}l}\|u_p\|_l
\sum_{-1\leq q\leq p}\lambda_q^{sl+\alpha-\frac{sl+\alpha}l}\|\theta_q\|_l^{l-1}\lambda_{q-p}^{1-\alpha+\frac{sl+\alpha}l}\\
\lesssim &c_0\kappa l^2\sum_{p\geq -1}\lambda_p^{sl+\alpha}\|u_p\|_l^l
\end{split}
\end{equation}
since $1-\alpha+\frac{sl+\alpha}l>0$;
\begin{equation}\notag
\begin{split}
D\lesssim &l^2f(t)\sum_{-1\leq p\leq Q}\lambda_p^{-1}\|u_p\|_l\sum_{-1\leq q\leq p}\lambda_q^{sl+1}\|\theta_q\|_l^{l-1}\\
\lesssim &l^2f(t)\sum_{-1\leq p\leq Q}\lambda_p^{s}\|u_p\|_l\sum_{-1\leq q\leq p}\lambda_q^{s(l-1)}\|\theta_q\|_l^{l-1}\lambda_{q-p}^{1+s}\\
\lesssim &l^2f(t)\sum_{-1\leq p\leq Q}\lambda_p^{sl}\|u_p\|_l^l
\end{split}
\end{equation}
for arbitrary $s>0$. Thus, we have
\begin{equation}\label{est-i3}
|I_3|\lesssim c_0\kappa l^2\sum_{p\geq -1}\lambda_p^{sl+\alpha}\|u_p\|_l^l+
l^2f(t)\sum_{-1\leq p\leq Q}\lambda_p^{sl}\|u_p\|_l^l.
\end{equation}
Combining (\ref{ineq-q})-(\ref{est-i3}) gives that 
\begin{equation}\notag
\frac{d}{dt}\sum_{q\geq -1}\lambda_q^{sl}\|\theta_q\|_l^l
\leq -C\kappa (1-c_0l^2)\sum_{q\geq -1}\lambda_q^{sl}\|\Lambda^{\frac 12} \theta_q\|_2^2+l^2f(t)\sum_{-1\leq  q\leq Q}\lambda_q^{sl}\|\theta_q\|_l^l.
\end{equation}
Therefore, one can choose small enough $c_0$ such that $c_0l^2<1/2$. It then follows from the Gronwall's inequality that 
\begin{equation}\notag
\|\theta(t)\|_{B^s_{l,l}}\leq l^2\|\theta(0)\|_{B^s_{l,l}}\exp\left({\int_0^tf(s)\,ds}\right).
\end{equation}
Thus $\theta$ has a finite norm in $B^s_{l,l}$ since $f\in L^1(0,T)$. It completes the proof of Theorem \ref{thm}.
%Notice that, to guarantee the embedding $W^{s,l}\subset C^{0, s-\frac2l}$ and hence claim the regularity of $\theta$ up to time $T$, it requires that $s-\frac 2l>1-\alpha$. 

\bigskip

Recall that in \cite{DP}, the authors obtained the regularity for the supercritical SQG provided $\theta\in L^s((0,T); B^\gamma_{p,\infty})$ with $\gamma=\frac2p+1-\alpha+\frac\alpha s$, for $1\leq s<\infty$ and $2< p<\infty$.

In the following lemma, we show that the regularity condition $f\in L^1(0,T)$ is weaker than the Prodi-Serrin criterion from \cite{DP}. It can be shown similarly that this condition is weaker than all the other Prodi-Serrin criteria mentioned in the introduction.

\begin{Lemma}\label{compare}
Let $\theta(t)$ be a weak solution to (\ref{QG}) on $[0,T]$. If $\theta\in L^s((0,T); B^\gamma_{p,\infty})$ with $\gamma=\frac2p+1-\alpha+\frac\alpha s$, for $1\leq s<\infty$ and $2< p\leq\infty$, then $f\in L^1(0,T)$.
\end{Lemma}
\pf 
Let $U=[0,T]\cap \{t:\Lambda(t)>1\}$. One can see that, by the definition of $\Lambda(t)$ in (\ref{wave}) and the Bernstein's inequality
\[\int_{[0,T]\textbackslash U}f(t) \, dt \lesssim \int_{[0,T]\textbackslash U}\|\theta(t)\|_2 \, dt<\infty.\]
Using (\ref{wave}) and the Bernstein's inequality again, we have
\begin{equation}\notag
\begin{split}
\int_Uf(t)\, dt=&\int_U\sup_{q\leq Q}\lambda_q\|\theta(t)\|_\infty\, dt\\
\leq& \int_U\Lambda(t)^{\alpha-\frac\alpha s}\sup_{q\leq Q}\lambda_q^{1-\alpha+\frac\alpha s}\|\theta(t)\|_\infty\, dt\\
\leq& \int_U\Lambda(t)^{\alpha-\frac\alpha s}\sup_{q\leq Q}\lambda_q^{\frac2p+1-\alpha+\frac\alpha s}\|\theta(t)\|_p\, dt\\
\leq&\left(\int_U\Lambda^\alpha(t)\, dt\right)^{1-\frac1s}\|\theta\|_{L^s(0,T;B^\gamma_{p,\infty})}\\
\leq &c(\kappa)\left(\int_U\Lambda^{\alpha+s-s\alpha}(t)\|\theta_{Q(t)}\|_p^s\, dt\right)^{1-\frac1s}\|\theta\|_{L^s(0,T;B^\gamma_{p,\infty})}\\
\leq &c(\kappa)\left(\int_U\Lambda^{\alpha+\frac{2s}p+s-s\alpha}(t)\|\theta_{Q(t)}\|_p^s\, dt\right)^{1-\frac1s}\|\theta\|_{L^s(0,T;B^\gamma_{p,\infty})}\\
\leq &c(\kappa)\|\theta\|_{L^s(0,T;B^\gamma_{p,\infty})}^s,
\end{split}
\end{equation}
which finishes the proof.

\cbdu

\begin{Remark}
In \cite{DP}, the authors were not able to obtain the regularity in the case of $p=\infty$, that is when $\theta\in L^s((0,T); B^\gamma_{\infty,\infty})$ with $\gamma=1-\alpha+\frac\alpha s$, for $1\leq s<\infty$. Lemma \ref{compare} and Theorem \ref{thm} imply that the regularity can be obtained provided $\theta$ is in such spaces.
\end{Remark}

\bigskip

\section{Energy conservation}
\label{sec-on}

In this section we prove Theorem \ref{thm-energy}. Formally, multiplying equation (\ref{QG}) by $\theta$ and integrating over space and time yields the energy equality
\begin{equation}\notag
\int_{\mathbb R^2}\frac12|\theta(t)|^2\, dx+\kappa\int_0^t\int_{\mathbb R^2}|\Lambda^{\alpha/2}\theta|^2\, dxd\tau=\int_{\mathbb R^2}\frac12|\theta(0)|^2\, dx.
\end{equation}
To prove a weak solution $\theta$ satisfies the energy equality on $[0,T]$, it is enough to show that the energy flux vanishes,
\[
\limsup_{Q \to \infty} \int_0^T\int_{\mathbb R^2}u\th \cdot\nabla (\theta_{\leq Q})_{\leq Q} \, dxdt=0,
\]
%provided $\theta\in L_{loc}^2(0,T;B^{1/2}_{2,c(\mathbb N)})$.
which can be shown similarly as in \cite{CDatt} (Section 3). For completeness, we present a brief proof in the following by using Littlewood-Paley decomposition method.
% Multiplying the first equation in (\ref{QG}) by $(\theta_{\leq Q})_{\leq Q}$ and integrating yields
%\begin{equation}\notag
%\frac{1}{2}\frac{d}{dt}\|\theta_{\leq Q}\|_{2}^2dx+\kappa\|\Lambda^{\frac 12}\theta_{\leq Q}\|_2^2=-\int_{\mathbb T^2}(u\cdot\nabla \theta)%\cdot (\theta_{\leq Q})_{\leq Q}dx+
%\int_{\mathbb T^2}f\cdot(\theta_{\leq Q})_{\leq Q}dx.
%\end{equation}
Denote 
\[r_Q(u,\theta)=\int_{\mathbb R^2}h_Q(y)\left(u(x-y)-u(x)\right) \left(\theta(x-y)-\theta(x)\right)dy.\]
Then 
\[(u\theta)_{\leq Q}=r_Q(u,\theta)-u_{>Q}\theta_{>Q}+u_{\leq Q}\theta_{\leq Q}.\]
The energy flux can be written as 
\begin{equation}\notag
\begin{split}
\Pi_Q=&\int_{\mathbb R^2}(u\theta)\cdot \nabla(\theta_{\leq Q})_{\leq Q}\,dx\\
=&\int_{\mathbb R^2}(u\theta)_{\leq Q}\cdot\nabla\theta_{\leq Q}\,dx\\
=&\int_{\mathbb R^2}r_Q(u,\theta)\cdot \nabla\theta_{\leq Q}dx
-\int_{\mathbb R^2}u_{>Q}\theta_{>Q}\cdot \nabla\theta_{\leq Q}\,dx.
\end{split}
\end{equation}
Here we used the fact that $\int_{\mathbb R^2}u_{\leq Q}\theta_{\leq Q}\cdot \nabla(\theta_{\leq Q})\,dx=0$ since $u$ is divergence free.
Since $\theta\in L^\infty(t_0,\infty;L^\infty)$ for all $t_0>0$ (see \cite{CV}, or \cite{CD}), there exists a constant $C$ such that $\|\theta(t)\|_\infty\leq C$ for all $t>t_0$. Thus,  combining the estimate
\begin{equation}\notag
\|u(\cdot-y)-u(\cdot)\|_2\lesssim\sum_{p\leq Q}|y|\lambda_p\|u_p\|_2+\sum_{p> Q}\|u_p\|_2,
\end{equation}
we have
\begin{equation}\notag
\begin{split}
\|r_Q(u,\theta)\|_2&\leq \int_{\mathbb R^2}h_Q(y)\|u(\cdot-y)-u(\cdot)\|_2\|\theta(\cdot-y)-\theta(\cdot)\|_\infty\,dy\\
&\lesssim \int_{\mathbb R^2}h_Q(y)\left(\sum_{p\leq Q}|y|\lambda_p\|u_p\|_2
+\sum_{p> Q}\|u_p\|_2\right)dy\\
&\lesssim \sum_{p\leq Q}\lambda_Q^{-1}\lambda_p\|u_p\|_2
+\sum_{p> Q}\|u_p\|_2.
\end{split}
\end{equation}
Therefore, applying Bernstein's inequality we obtain that
\begin{equation}\notag
\begin{split}
\left|\int_{\mathbb R^2}r_Q(u,\theta)\cdot\nabla\theta_{\leq Q}\,dx\right|
\leq& \|r_Q(u,\theta)\|_2\|\nabla \theta_{\leq Q}\|_2\\
\lesssim& \sum_{p\leq Q}\lambda_Q^{-1}\lambda_p\|u_p\|_2\sum_{p'\leq Q}\lambda_{p'}\|\theta_{p'}\|_2\\
&+\sum_{p> Q}\|u_p\|_2\sum_{p'\leq Q}\lambda_{p'}\|\theta_{p'}\|_2\\
:= &I+II.
\end{split}
\end{equation}
%Since $\theta\in L^\infty(t_0,\infty;L^\infty)$ for all $t_0>0$ (see \cite{CV}, or \cite{CD}), it implies $u\in L^{\infty}(t_0,\infty;BMO)$. Therefore there exists $C$, such that $\|u_q(t)\|_\infty \leq C$ for all $q$ and a.a. $t >t_0$.
Due to the fact $\|u_p\|_2\lesssim \|\theta_p\|_2$ for all $p\geq -1$,  $I$ and $II$ are estimated as, by using H\"older's inequality, Young's inequality and Jensen's inequality
\begin{equation}\notag
\begin{split}
I\lesssim &\sum_{p\leq Q}\lambda_Q^{-1}\lambda_p\|\theta_p\|_2\sum_{p'\leq Q}\lambda_{p'}\|\theta_{p'}\|_2\\
\lesssim &\sum_{p\leq Q}\lambda_{p-Q}^{\frac 12}\lambda_p^{\frac 12}\|\theta_p\|_2\sum_{p'\leq Q}\lambda_{p'-Q}^{\frac 12}\lambda_{p'}^{\frac 12}\|\theta_{p'}\|_2\\
\lesssim &\left(\sum_{p\leq Q}\lambda_{p-Q}^{\frac12}\lambda_p^{\frac12}\|\theta_p\|_2\right)^2+\left(\sum_{p'\leq Q}\lambda_{p'-Q}^{\frac 12}\lambda_{p'}^{\frac 12}\|\theta_{p'}\|_2\right)^2\\
\lesssim &\sum_{p\leq Q}\lambda_{p-Q}^{\frac12}\lambda_p\|\theta_p\|_2^2+\sum_{p'\leq Q}\lambda_{p'-Q}^{\frac 12}\lambda_{p'}\|\theta_{p'}\|_2^2,
\end{split}
\end{equation}
and
\begin{equation}\notag
\begin{split}
II\lesssim &\sum_{p> Q}\|\theta_p\|_2\sum_{p'\leq Q}\lambda_{p'}^2\|\theta_{p'}\|_2\\
\lesssim &\sum_{p> Q}\lambda_{Q-p}^{\frac 12}\lambda_p^{\frac 12}\|\theta_p\|_2\sum_{p'\leq Q}\lambda_{p'-Q}^{\frac 12}\lambda_{p'}^{\frac 12}\|\theta_{p'}\|_2\\
\lesssim &\left(\sum_{p> Q}\lambda_{Q-p}^{\frac12}\lambda_p^{\frac12}\|\theta_p\|_2\right)^2
+\left(\sum_{p'\leq Q}\lambda_{p'-Q}^{\frac 12}\lambda_{p'}^{\frac 12}\|\theta_{p'}\|_2\right)^2\\
\lesssim &\sum_{p> Q}\lambda_{Q-p}^{\frac12}\lambda_p\|\theta_p\|_2^2
+\sum_{p'\leq Q}\lambda_{p'-Q}^{\frac 12}\lambda_{p'}\|\theta_{p'}\|_2^2.
\end{split}
\end{equation}
On the other hand, we have the similar estimate
\begin{equation}\notag
\begin{split}
\left|\int_{\mathbb R^2}u_{>Q}\theta_{>Q}\cdot\nabla\theta_{\leq Q}\,dx\right|
\lesssim &\sum_{p>Q}\|\theta_p\|_2\sum_{p'\leq Q}\lambda_{p'}\|\theta_{p'}\|_2\\
\lesssim &\sum_{p>Q}\lambda_{Q-p}^{\frac 12}\lambda_{p}^{\frac 12}\|\theta_p\|_2\sum_{p'\leq Q}\lambda_{p'-Q}^{\frac 12}\lambda_{p'}^{\frac 12}\|\theta_{p'}\|_2\\
\lesssim &\left(\sum_{p> Q}\lambda_{Q-p}^{\frac12}\lambda_p^{\frac12}\|\theta_p\|_2\right)^2
+\left(\sum_{p'\leq Q}\lambda_{p'-Q}^{\frac12}\lambda_{p'}^{\frac12}\|\theta_{p'}\|_2\right)^2\\
\lesssim &\sum_{p> Q}\lambda_{Q-p}^{\frac12}\lambda_p\|\theta_p\|_2^2
+\sum_{p'\leq Q}\lambda_{p'-Q}^{\frac12}\lambda_{p'}\|\theta_{p'}\|_2^2.
\end{split}
\end{equation}
Therefore 
\begin{equation}\notag%\label{flux}
\left|\Pi_Q\right|\lesssim \sum_{p\leq Q}\lambda_{p-Q}^{\frac12}\lambda_p\|\theta_p\|_2^2
+\sum_{p> Q}\lambda_{Q-p}^{\frac12}\lambda_{p}\|\theta_{p}\|_2^2
= \sum_{p\geq -1}\lambda_{|p-Q|}^{-{\frac12}}\lambda_p\|\theta_p\|_2^2.
\end{equation}
Since $\theta\in L^2(0,T;B^{1/2}_{2,c(\mathbb N)})$, we have 
\begin{equation}\notag
\lim\sup_{Q\to\infty}\int_0^T|\Pi_Q|\, dt\lesssim \lim\sup_{Q\to\infty}\int_0^T\lambda_Q\|\theta_Q\|_2^2\, dt=0.
\end{equation}
It completes the proof of Theorem \ref{thm-energy}.
%Suppose $\theta\in L_{loc}^2(0,T;B^{1/2}_{2,c(\mathbb N)})$, we have that
%\[\int_{t_0}^{T}\lambda_p\|\theta_p\|_2^2\,dt\to 0 \qquad \mbox { as } p\to\infty.\]
%It implies the right hand side of (\ref{flux}) converges to 0 as $Q\to\infty$, which gives that 
%\[ \lim_{Q\to \infty}\int_{t_0}^{T} \Pi_Q \,dt=0. \]

%%%%%%%%%%%%%%%%%%%%%%%%%%%%%%%

 %Proof of the last inequality: $\lim\sup_{Q\to\infty}\int_0^T|\Pi_Q|\, dt\lesssim \lim\sup_{Q\to\infty}\int_0^T\lambda_Q\|\theta_Q\|_2^2\, dt$: 
 
 %Let $a_q=\lambda_q\|\theta_q\|_2^2$. Let $M=\lim\sup_{Q\to\infty} a_Q$. For any $\epsilon>0$, there exists $Q_0$ such that $a_q>M-\epsilon$ for $q>Q_0$. Then for large enough $Q^*$ such that $Q_0\ll Q^*$, we have
%\begin{equation}\notag
%|\Pi_{Q^*}|\lesssim \sum_{p\leq Q_0}\lambda_{|p-Q^*|}^{-1}\lambda_p\|\theta_p\|_2^2+\sum_{p\geq Q_0}\lambda_{|p-Q^*|}^{-1}\lambda_p\|\theta_p\|_2^2\lesssim 2M+M\sum_{p\geq Q_0}\lambda_{|p-Q^*|}^{-1}\lambda_p\lesssim M.
%\end{equation}

%\Endrefs
\end{document}